\documentclass[12pt
]{amsart}
\usepackage{amssymb}
\usepackage[dvips]{graphics}
\ifx\mathrm\undefined\let\mathrm\rm\fi
\ifx\mathbf\undefined\let\mathbf\bf\fi
\ifx\mathfrak\undefined\let\mathfrak\frak\fi
\ifx\mathcal\undefined\let\mathcal\cal\fi
\ifx\mathbb\undefined\let\mathbb\Bbb\fi
\ifx\emph\undefined\let\emph\it\fi
 at9.98pt

\newcommand{\SL}{\mathrm{SL}}

\newcommand{\Z}{{\mathbb Z}}

\newcommand{\R}{{\mathbb R}}
\newcommand{\C}{{\mathbb C}}

\newcommand{\Ref}[1]{{(\ref{#1})}}

\newcommand{\sig}{{\sigma}}

\newcommand{\dontprint}[1]
{\relax}

\newtheorem%
{thm}{Theorem}[section]
\newtheorem%
{proposition}[thm]{Proposition}
\newtheorem%
{lemma}[thm]{Lemma}
\newtheorem%
{lemmadef}[thm]{Lemma-Definition}
\newtheorem%
{corollary}[thm]{Corollary}
\newtheorem%
{conjecture}[thm]{Conjecture}

\newcommand{\bea}{\begin{eqnarray*}}
\newcommand{\eea}{\end{eqnarray*}}
\newcommand{\bean}{\begin{eqnarray}}
\newcommand{\eean}{\end{eqnarray}}

\newcommand{\mdpd}[4]{
\left(\begin{array}{rr}#1 & #2 \\ #3 & #4\end{array}\right)}

\newcommand{\HH}{{\mathcal{H}}}

\begin{document}
\title[{}]{Hypergeometric theta functions and elliptic Macdonald polynomials}

\author[{}]
{Giovanni Felder${}^{*,1}$ 
\and Alexander Varchenko${}^{**,2}$}
\thanks{${}^1$Supported in part by the Swiss National Science Foundation} 
\thanks{${}^2$Supported in part by NSF grant DMS-0244579}
\maketitle 
\medskip \centerline{\it ${}^*$Department of Mathematics,
  ETH-Zentrum,} \centerline{\it 8092 Z\"urich, Switzerland} \medskip
\centerline{\it ${}^{**}$Department of Mathematics, University of
  North Carolina at Chapel Hill,} \centerline{\it Chapel Hill, NC
  27599-3250, USA} \medskip 

\centerline{August, 2003}
\begin{abstract}
Elliptic Macdonald polynomials of $sl_2$-type and level 2 are introduced.
Suitable limits of elliptic Macdonald polynomials are the standard Macdonald 
polynomials and conformal blocks.
 Identities for elliptic Macdonald polynomials, in particular their modular properties,
 are studied.

\end{abstract}

\section{Introduction}

Conformal field theory and its q-deformation provide new classes of
interesting special functions. Here we study special functions arising in
the q-deformation of conformal blocks on elliptic curves with one
marked point. These functions are given by 
elliptic versions of q-hypergeometric integrals and we call them
hypergeometric theta functions.
We study here the simplest non-trivial case, corresponding to
the 3-dimensional irreducible representation of the Lie algebra $sl_2$,
but one should expect the same picture to hold in the more general
situation of the $m$-th symmetric tensor power of the vector representation
of $sl_{n+1}$. Then hypergeometric theta functions are functions
on the Cartan subalgebra of $sl_{n+1}$ depending on an integer
parameter $\kappa$, the level, and two complex parameters:
$\tau$, parametrizing the elliptic curve and $\eta$, the
deformation parameter. In the {\em classical limit} $\eta\to 0$, hypergeometric
theta functions are supposed to 
converge to conformal blocks, realized as theta functions
obeying vanishing conditions and obeying the KZB equation. In the
{\em trigonometric limit} $\tau\to i\infty$ one expects to
recover Macdonald polynomials.
More precisely, the quotient of hypergeometric theta functions by
a suitable product of ordinary theta functions is expected to 
converge to $A_n$-Macdonald polynomials in the trigonometric limit.
We call these quotients {\em elliptic Macdonald polynomials}. From the point
of view of representaton theory, hypergeometric theta functions are
expected to be related to traces of intertwining operators for
the quantum affine Lie algebra $U_q(\hat {sl}_{n+1})$, $q=e^{-2\pi i\eta}$,
over integrable simple modules of
level $\kappa-n-1$. Thus we conjecture that
our definition of elliptic Macdonald polynomials
coincides with the definition of Etingof and
Kirillov \cite{EK2} in terms of traces of intertwining operators. 

The classical limit is the limit in which classical Lie algebras are
recovered from their quantum version, and in the trigonometric limit
one goes from affine to simple finite dimensional Lie algebras. Thus, we
expect the following picture to hold.
   \[
\begin{array}{ccc}
\text{Conformal blocks}&
\stackrel
{\tau\to i\infty}{\longrightarrow}&
\text{Jack polynomials}\\
\uparrow \eta\to 0
&&\uparrow \eta\to 0
\\
\text{Hypergeometric theta functions}
&\stackrel
{\tau\to i\infty}{\longrightarrow}& \text{$A_n$-Macdonald polynomials}
\end{array}
\]
In this paper, we establish this picture in the simplest non-trivial
$sl_2$-case. For $sl_2$, the
hypergeometric theta functions may  be considered as
a degeneration of the elliptic hypergeometric integrals
studied for generic parameters in \cite{FV2}, \cite{FV3}.
The elliptic hypergeometric integrals obey several identities
involving their values at points related to each other by an
action of $\mathrm{SL}(3,\mathbb Z)$. Some of these identities
survive in the degenerate case and result in the known identities
in the trigonometric and classical limits.
In the case we consider, these identities are the following:
first of all, we see that the KZB equation obeyed by
conformal blocks is a limiting case of two equations, the qKZB equations,
 both of which are solved
by hypergeometric theta functions. One is an integral equation,
which in the trigonometric limit becomes a Macdonald-Mehta type identity
for Macdonald polynomials. The other is an infinite difference equation,
which in the trigonometric limit expresses the fact that Macdonald
polynomials are eigenvectors of the Macdonald--Ruijsenaars difference
operators.
Other identities in \cite{FV3} 
give the transformation properties of hypergeometric theta
functions under the modular group $\mathrm{SL}(2,\Z)$.
In the classical limit, we recover the known modular
properties of conformal blocks.
Finally, the orthogonality relations of \cite{FV3} also have
a counterpart for elliptic Macdonald polynomials and reduce
to the characterizing orthogonality properties of Macdonald
polynomials.

\bigskip

\section{qKZB heat equation}

\subsection{Theta functions of level $\kappa$}
Let $\mathrm{Im\,\tau}>0$.
A holomorphic function $f:\C\to \C$ is called
a theta function of level $\kappa\in\Z_{\geq0}$ if
\[
f(\lambda+2r+2s\tau)=e^{-2\pi i\kappa(s^2\tau+s\lambda)}
f(\lambda),\qquad r,s\in\Z.
\]
Let $\Theta_\kappa(\tau)$ denote the space of
theta functions of level $\kappa$. If $\kappa=0$,
$\Theta_\kappa(\tau)=\C$. If $\kappa>0$,
$\Theta_\kappa(\tau)$  has dimension
$2\kappa$: a basis is
\[
\theta_{j,\kappa}(\lambda,\tau)=\sum_{n\in\Z+\frac j{2\kappa}}
e^{2\pi i\kappa (n^2\tau+n\lambda)},\qquad j\in\Z/2\kappa\Z.
\]
The space $\Theta_\kappa(\tau)$ is the direct sum of
the space $\Theta_\kappa(\tau)^\mathrm{even}$
of even theta functions and the
space $\Theta_\kappa(\tau)^\mathrm{odd}$
of odd theta functions, of
dimension $\kappa+1$ and $\kappa-1$, respectively. In
particular, for $\kappa=2$, Jacobi's first
theta function
$
\theta(\lambda,\tau)=-\sum_{n\in\Z+1/2}
e^{\pi i(n^2\tau+n(2\lambda+1))}
$,
($\theta=i(\theta_{-1,4}-\theta_{1,4})$)
spans the space of odd theta functions of level 2.

\bigskip

\subsection{The elliptic hypergeometric integral}
Let $\lambda,\mu,\tau,\sigma, \eta\in\C$ and
$\mathrm{Im}\,\tau >0$, $\mathrm{Im}\,\sigma >0$, 

The basic object is the hypergeometric integral
\[u(\lambda,\mu,\tau,\sig,\eta)=
e^{-\frac{\pi i\lambda\mu}{2\eta}}\int_\gamma
\Omega_{2\eta}(t,\tau,\sig)
\frac{\theta(\lambda+t,\tau)\theta(\mu+t,\sig)}
     {\theta(t-2\eta,\tau)\theta(t-2\eta,\sig)}dt.
\]
where
\[
\Omega_{2\eta}(t,\tau,\sig)=\prod_{j,k=0}^\infty
\frac
{(1-e^{2\pi i(t-2\eta +j\tau+k\sig)})(1-e^{2\pi i(-t-2\eta +(j+1)\tau+(k+1)\sig)})}
{(1-e^{2\pi i(t+2\eta +j\tau+k\sig)})(1-e^{2\pi i(-t+2\eta +(j+1)\tau+(k+1)\sig)})}
\,.
\]
The integration path is the interval $[0,1]$ as long as
$\eta$ has positive imaginary part. For general $\eta$
the function $u$ is defined by analytic continuation. The analytic
properties of $u$ are given by the following refinement
of a result in \cite{FV3}.

\bigskip

\begin{proposition}\label{p-re}
Let $H_+=\{z\in\C\,|\,\mathrm{Im}\,z>0\}$  be the
upper half plane.
Then $u$ is a meromorphic function on $\C\times\C\times 
H_+\times H_+\times\C$. It is regular on the complement
of the hyperplanes $4\eta+l+m\tau+n \sig=0$,
 $l\in\Z, m,n\in\Z_{\geq0}$,
where it has at most simple poles.
\end{proposition}

\bigskip

\noindent{\it Proof:} The integrand, regarded as a
function of $t$ for fixed generic $\tau,\sig$,
has two families of simple poles: the first family 
consists of the poles at 
$t=2\eta+l+m\tau+n \sig$ and the second at
$t=-(2\eta+l+m\tau+n\sig)$. In both families $l,m,n$
run over integers such that $m,n \geq 0$.
If $\eta$ has positive imaginary part, the poles
of the first family lie in the upper half plane and
the poles of the second family lie in the lower half plane. 
The integration cycle separates the (projection onto 
$\C/\Z$ of) the two families.
As we do an analytic continuation, the integration
cycle, originally along the real axis, gets deformed
and $u$ is holomorphic as long as no pole of the first
family coincides with a pole of the second family.
Such a coincidence happens when $4\eta+l+m\tau+n \sig=0$
for some integers $l,m,n$ with $m,n\geq0$. As
one approaches a generic point of a hyperplane
$4\eta+l+m\tau+n\sig=0$, a finite number of poles of
the first family meets poles of the second, and the
integration cycles gets pinched between simple poles.
The integral has then at most a simple pole on the
hyperplane. \hfill$\square$

We shall use the following two properties of the hypergeometric integral.

\bigskip

\begin{lemma}\label{sym of u}
The hypergeometric integral satisfies the relation
\bea
u(-\lambda,-\mu,\tau,\sig,\eta)\ =\ u(\lambda,\mu,\tau,\sig,\eta)\ .
\eea
\end{lemma}

\bigskip

The lemma follows from identity (A.15) in \cite{FV3}, cf. identities in \cite{FV4}.

\bigskip

\begin{lemma}[\cite{FV3}]\label{vanishing} For $r,s\in \Z$ we have
\bea
u(\lambda, 2\eta + r + sp, \tau,\sig,\eta)\ =\ 
e^{2\pi i s (\tau - 4\eta)}\
u(\lambda,-2\eta + r + sp,\tau,\sig,\eta)\ .
\eea
\end{lemma}

\bigskip

\subsection{The qKZB discrete connection}
Let 
\[
Q(\mu,\sig,\eta) =
\frac{\theta(4\eta,\sig)\theta'(0,\sig)}
{\theta(\mu-2\eta,\sig)
  \theta(\mu+2\eta,\sig)}\,,
\]
where the prime denotes the derivative with respect to
the first argument.
In \cite{FV2}, \cite{FV3} we studied the integral
operator (defined on a suitable space of holomorphic functions)
\begin{equation}\label{e-U}
U(\tau,\sig,\eta):v\mapsto \int_{\eta\R} u(\lambda,\mu,\tau,\sig,\eta)
Q(\mu,\sig,\eta)v(-\mu)d\mu,
\end{equation}
for generic values
of the parameters. In \cite{FV2} this integral
operator was used to define a q-deformation of
the KZB heat equation, the qKZB heat equation. 
In \cite{FV3} the integral operator was used
to construct solutions of the qKZB heat equations
and describe their monodromy.

In this paper  we consider the case when $\sigma$ is a positive integer
multiple of $-2\eta$. In this case the qKZB heat operator
can be defined as a map between finite dimensional
vector spaces of theta functions.

For $\kappa\in\Z$, $\kappa\geq4$, let $E_\kappa(\tau,\eta)$
be the subspace of $\Theta_{\kappa+2}(\tau)^\mathrm{odd}$
consisting of functions vanishing at $2\eta j+\Z+\Z\tau$,
$j=-1,0,1$.
Since such functions are divisible
by the odd function
$\theta(\lambda-2\eta,\tau)\theta(\lambda,\tau)
\theta(\lambda+2\eta,\tau)$ in the 
graded 
ring $\oplus_{\kappa=0}^\infty \Theta_\kappa(\tau)$, we have
 \begin{equation}\label{e-cb}
E_\kappa(\tau,\eta)=\{
\theta(\lambda-2\eta,\tau)\theta(\lambda,\tau)
\theta(\lambda+2\eta,\tau)g(\lambda)\ |\
g\in\Theta_{\kappa-4}(\tau)^\mathrm{even}\}
\end{equation}
In particular,
\[
\mathrm{dim}\,E_\kappa(\tau,\eta)=\kappa-3.
\]

\bigskip

\begin{proposition} \label{map of thetas}
Assume that
$\kappa$ is an integer $\geq 4$,
$\mathrm{Im}\,\tau>0,\mathrm{Im}\,\eta<0$,
$j\tau+4\eta\not\in\Z$, for every $ j = 1, 2, \dots$.
Let $\alpha(\lambda,\eta)=\exp(-{\pi i\lambda^2/4\eta})$.
Denote by $\alpha(\eta)$ the operator of multiplication
by $\alpha(\lambda,\eta)$.
Then the integral operator
\[
T_\kappa(\tau,\eta)=-\frac{e^{4\pi i\eta}}
{2\pi\sqrt{4i\eta}}\,
\alpha(\eta)\, U(\tau,\tau-2\eta\kappa,\eta)\,
\alpha(\eta)\]
 maps 
$E_\kappa(\tau-2\eta\kappa,\eta)$ to $E_\kappa(\tau,\eta)$.
\end{proposition}

\bigskip

\begin{proof}
Our goal is to show that for any function $f \in E_\kappa(\tau - 2\eta \kappa)$
the integral
\bea
-\frac{e^{4\pi i\eta}}{2\pi\sqrt{4i\eta}}\,
\int_{\eta \R}\int_\gamma
e^{-\pi i\frac{(\lambda+\mu)^2}{4\eta}}
\Omega_{2\eta}(t,\tau,\tau-2\eta \kappa)
\frac{\theta(\lambda+t,\tau)\theta(\mu+t,\tau-2\eta \kappa)}
     {\theta(t-2\eta,\tau)\theta(t-2\eta,\tau-2\eta \kappa)}\,
\\
\ {}\ \ {}\ \ {}\ \ {}\ \ {}\ \ {}\ \ {}\ \ {}\ \ {}\ \ {}\ \ {}\ \ {}\ \times
\frac{\theta(4\eta,\tau-2\eta \kappa)\theta'(0,\tau-2\eta \kappa)}
{\theta(\mu-2\eta,\tau-2\eta \kappa)
  \theta(\mu+2\eta,\tau-2\eta \kappa)}\,f(-\mu)\,
dt d\mu
\eea
converges and belongs to $E_\kappa(\tau )$ as a function of $\lambda$.

Clearly the integrand is holomorphic in $\mu$.
The integrand has the same 
poles with respect to $t$ as the integrand of 
$u(\lambda, \mu ,\tau, \tau -2\eta\kappa,\eta)$, so
the integration with respect to $t$ is defined  except possibly
on the hyperplanes of Proposition \ref{p-re}. Under
the assumptions of this proposition that condition 
reduces to the condition $4\eta+j\tau\not\in\Z$, \, $j=1,2,\dots$.

\begin{lemma}[cf. Lemma C.1, \cite{FV3}]
For $f\in \Theta_\kappa(\tau)$, there exists $C_1, C_2 > 0$ such that
for all $\lambda \in \C$, we have
\bea
|f(\lambda)| \ \leq\ C_1\,\exp \left(\, \frac{\pi \kappa}2\,
\frac {(\mathrm{Im}\, \lambda)^2}{\mathrm{Im}\, \tau}\,+\,
 C_2 \,\mathrm{Im}\, \lambda\,\right)\ .
\eea
\end{lemma}
It follows from the lemma that if $\mu = \eta x,\, x\in \R$, then
there exists $C_3, C_4 > 0$ such that the absolute value of the integrand 
is not greater than
\bea
C_3\, \exp \left(\,\frac{\pi} 4\, x^2\, ( \mathrm{Im}\,\eta )\,
\frac {\mathrm{Im}\,\tau - 4\,\mathrm{Im}\,\eta}{\mathrm{Im}\,\tau- 2\kappa\,\mathrm{Im}\,\eta}\,+\,
 C_4 \, |x|\,\right)\ .
\eea
Thus the integration with respect to $\mu$ is defined and the double integral
is a holomorphic function of $\lambda$.

The integral is an odd function of $\lambda$ by Lemma \ref{sym of u}.

The fact that the integral belongs to 
$\Theta_{\kappa+2}(\tau)$ easily follows from the theta function properties 
of the integrand. The fact that the integral is divisible by
$\theta(\lambda-2\eta,\tau)\theta(\lambda,\tau)
\theta(\lambda+2\eta,\tau)$ follows from Lemma \ref{vanishing}.

\end{proof}

We say that a function $v(\lambda,\tau)$ is a solution
of the qKZB heat equation if 
\begin{equation}\label{e-qKZB}
T_\kappa(\tau,\eta)v(\lambda,\tau-2\eta\kappa)=
v(\lambda,\tau),
\end{equation} 
Solutions $v(\lambda,\tau)$ belonging to
$E_\kappa(\tau,\eta)$ for each fixed $\tau$
are constructed below.
\section{Hypergeometric theta functions}

\subsection{Definition}
Let $\kappa\in\Z$, $\kappa\geq 4$. 
We say that an integer $l$ is admissible with respect to $\kappa$ if
 $l \neq \pm 1 \mod \kappa$.

For an admissible $l$, define the {\em $l$-th non-symmetric
hypergeometric theta function of level $\kappa+2$} by
\bea
\widetilde\Delta_{l,\kappa}(\lambda,\tau,\eta) \ = \
\sum_{j\in 2\kappa\Z+l}
\,u(\lambda,2\eta j,\tau,-2\eta\kappa,\eta)
 \ Q(2\eta j, -2\eta\kappa, \eta) \ 
e^{{\pi i}\frac{\tau + 4\eta}{2\kappa}j^2}\ .
\eea
Using the
transformation properties of $Q$ we may rewrite this as
\bea
\widetilde\Delta_{l,\kappa}(\lambda,\tau,\eta)  = e^{\frac{4\pi i \eta}{\kappa}l^2}
\ Q(2\eta l,-2\eta\kappa,\eta) \sum_{j\in 2\kappa\Z+l}
u(\lambda,2\eta j,\tau,-2\eta\kappa,\eta)
\ e^{{\pi i}\frac{\tau - 4\eta}{2\kappa}j^2}\ .
\eea
We have
\bea
\widetilde\Delta_{l,\kappa}(-\lambda,\tau,\eta)  = 
\widetilde\Delta_{-l,\kappa}(\lambda,\tau,\eta) \ , 
\qquad
\widetilde\Delta_{l,\kappa}(\lambda,\tau,\eta)  = 
\widetilde\Delta_{l+2\kappa,\kappa}(\lambda,\tau,\eta)\ .
\eea
We define the {\em $l$-th hypergeometric theta function of level $\kappa+2$} by
\bea
\Delta_{l,\kappa}(\lambda,\tau,\eta)  = \widetilde\Delta_{l,\kappa}(\lambda,\tau,\eta) -
\widetilde\Delta_{l,\kappa}(-\lambda,\tau,\eta)\ .
\eea

\bigskip

\begin{thm}\label{t-1}  Let $\mathrm{Im}\,\eta<0$. Then        \ %
\begin{enumerate} 
\item[(i)] For any fixed $\tau\in H_+$, such
that $j\tau+4\eta\not\in\Z$ for every
 $j=1,2,\dots$, and for any admissible $l$,
the series $\widetilde\Delta_{l,\kappa}$
converges to a holomorphic function 
$\lambda\mapsto\widetilde\Delta_{l,\kappa}(\lambda,\tau,\eta)$.
\item[(ii)] Under the same assumptions, we have
\bea\label{int repn}
\widetilde\Delta_{l,\kappa}(\lambda,\tau,\eta)\ =\
e^{\frac{2\pi i\eta }\kappa\, l^2} I_{l,\kappa}(\lambda,\tau,\eta)
Q(2\eta l,-2\eta\kappa,\eta)\ ,
\eea
where
\bea\label{int repn 2}
I_{l,\kappa}(\lambda,\tau,\eta)=\int_\gamma
\Omega_{2\eta}(t,\tau,-2\eta\kappa)
\frac{\theta(\lambda+t,\tau)\theta(2\eta l+t,-2\eta\kappa)}
     {\theta(t-2\eta,\tau)\theta(t-2\eta,-2\eta\kappa)}
e^{-2\pi ilt/\kappa}
\theta_{l,\kappa}\left({\textstyle\frac2\kappa}\,t-\lambda,\tau
\right)dt.
\eea
\item[(iii)] The functions  $\Delta_{l,\kappa}$,
$l=2,\dots,\kappa-2$, form a basis of $E_\kappa(\tau,\eta)$.
\end{enumerate}
\end{thm}

\bigskip 

\noindent{\it Proof:}
One shows that the integral $I_{l,\kappa}(\lambda,\tau,\eta)$
is defined and holomorphic with respect to $\lambda$
as in the proof of Proposition \ref{map of thetas}.
 The series for $\widetilde\Delta_{l,\kappa}$ is
then obtained by expanding the theta function $\theta_{l,\kappa}$
in a Fourier series in $\lambda$. This proves (i) and (ii).

The fact that 
 $I_{l,\kappa}(\lambda,\tau,\eta)$
(and thus $\widetilde\Delta_{l,\kappa}$) is a theta function of
level $\kappa+2$ follows easily from the theta function properties
of $\theta$, $\theta_{l,\kappa}$. The vanishing condition
for $\Delta_{l,\kappa}$ at the translates of $\pm 2\eta,0$ follows from 
Lemma \ref{vanishing}.
To complete the proof of (iii) it remains to show that 
the functions $\Delta_{l,\kappa}$,
$l=2,\dots,\kappa-2$, are linearly independent.
This follows from a degenerate version
of the {\em inversion relation} of \cite{FV3}:
\begin{equation}\label{e-ir}
\frac1{32\pi^2\eta}\,\int_0^2
u(-\mu,2\eta l,\tau,\sig,-\eta)\,u(\mu,2\eta n,\tau,\sig,\eta)
\,Q(\mu,\tau,\eta)\,d\mu
\ =\ Q(2\eta n,\sig,\eta)^{-1}\,\delta_{l,n}\ ,
\end{equation}
for any $l,n\in\Z$. The integration is along a path which
does not cross the straight line segment between $2\eta$
and $-2\eta$ or its translates by $\Z+\sig\Z$. 
 The proof of this
identity is the same as the proof of the inversion relation
in \cite{FV3}. The only difference is that the integration
is over a period instead of an infinite line, which is 
permitted since the integrand is a 2-periodic function
of $\mu$ for integers $l,n$.

We now set $\sig=-2\eta\kappa$ in the inversion relation
and restrict to
$n\not\equiv \pm1\mod 2\kappa$ to avoid poles of $Q$. 
We get, for $l, j=2,\dots,\kappa-2$,
\bean\label{inversion}
&&
\frac1{32\pi^2\eta}\,\int_0^2
\,u(-\mu,2\eta l,\tau,-2\eta\kappa,-\eta)\,
\Delta_{j,\kappa}(\mu,\tau,\eta)\,
Q(\mu,\tau,\eta)\,d\mu\
\ {}\ \ {}\ \ {}\ \ {}\ \ {}\ \ {}\ \ {}\ \ {}\ \ {}\ \ {}\ \ {}\ \ {}\ 
\\
&&
\ {}\ \ {}\ \ {}\ \ {}\ \ {}\ \ {}\ \ {}\ \ {}\ \ {}\ \ {}\ \ {}\ \ {}\ 
\ {}\ \ {}\ \ {}\ \ {}\ \ {}\ \ {}\ \ {}\ \ {}\ \ {}\ \ {}\ \ {}\ \ {}\ 
\ {}\ \ {}\ \ {}\ \ {}\ \ {}\ \ {}\ \ {}\ \ {}\ \ {}\ \ {}\ \ {}\ \ {}\ 
=\
\delta_{l,j}\,e^{\pi i\frac{4\eta+\tau}{2\kappa}j^2}\ .
\notag
\eean
In particular, this implies the linear independence of 
$\Delta_{l,\kappa}$.
\hfill$\square$

\subsection{Theta function solutions of the qKZB heat equation}

\begin{thm}\label{qkzb solutions}
For any admissible $l$ 
the functions $\widetilde\Delta_{l,\kappa}$ and
$\Delta_{l,\kappa}$ are solutions of the
qKZB heat equation:
\bea
T_\kappa(\tau,\eta)
\widetilde\Delta_{l,\kappa}(\lambda,\tau-2\eta\kappa,\eta)
&=&
\widetilde\Delta_{l,\kappa}(\lambda,\tau,\eta)\ ,
\\
T_\kappa(\tau,\eta)
\Delta_{l,\kappa}(\lambda,\tau-2\eta\kappa,\eta)
&=&
\Delta_{l,\kappa}(\lambda,\tau,\eta)\ .
\eea
\end{thm}

The theorem follows from the fact that for a fixed $\mu$
the function $u(\lambda,\mu,\tau-2\eta\kappa,-2\eta\kappa,\eta)$ gives
a solution of the qKZB equation up to a scalar factor depending on $\mu$,
see \cite{FV3}. Multiplying by the 
exponential function in the definition
of $\widetilde\Delta_{l,\kappa}$ we get rid of the scalar factor,
so that each term in the series obeys the qKZB equation.

\subsection{Theta functions as eigenfunctions of a difference operator}
The integral operator $U(\tau,\sigma,\eta)$ has a discrete version,
\bea
\bar U(\tau,\sigma,\eta) : v \mapsto
\sum_{m\in \Z}
u(\lambda, -\lambda+2\eta m,\tau, \sigma, \eta)
Q(-\lambda+2\eta m,\sigma,\eta)
v(\lambda-2\eta m) .
\eea

\begin{thm}

Assume that $\kappa$ is an integer $\geq 4$,
$\mathrm{Im}\,\tau>0,\,\mathrm{Im}\,\eta<0$,
$j\tau+4\eta\not\in\Z$ for every $ j = 1, 2, \dots$.
Then the difference operator
\bean\label{difference oper}
\bar T_\kappa(\tau,\eta) &=& 
C
\alpha(\eta)\, \bar U(\tau,\tau-2\eta\kappa,\eta)\,
\alpha(\eta)\ ,
\\
C &=&  \frac{ie^{4\pi i\eta}}{2\pi}\,
( \sum_{m\in \Z} e^{-\pi i \eta m^2})^{-1}\ ,
\notag
\eean
is well defined on $E_\kappa(\tau-2\eta\kappa,\eta)$
and  maps 
$E_\kappa(\tau-2\eta\kappa,\eta)$ to $E_\kappa(\tau,\eta)$.
Moreover, for any $l = 2, \dots , \kappa-2$ we have
\bean\label{difference}
\bar T_\kappa(\tau,\eta)
\widetilde \Delta_{l,\kappa}(\lambda,\tau-2\eta\kappa,\eta)
&=&
\widetilde \Delta_{l,\kappa}(\lambda,\tau,\eta) ,
\\
\bar T_\kappa(\tau,\eta)
 \Delta_{l,\kappa}(\lambda,\tau-2\eta\kappa,\eta)
&=&
 \Delta_{l,\kappa}(\lambda,\tau,\eta) .
\notag
\eean
\end{thm}

The proof is the same as the proof of Theorem 2.1 in \cite{FV3}, cf. Remark in
Section 3.3 in \cite{FV2}.

\section{Modular transformation properties}
\subsection{Formulas}

\begin{thm}\label{t-mod}
\begin{enumerate}
\item[(i)]
If \, $\mathrm{Im}\,\eta<0$, \, $\mathrm{Im}\,\tau>0$, then
for $l=2, \dots , \kappa-2$ we have
\bean\label{1mp}
\Delta_{l,\kappa}(\lambda,\tau+1,\eta)\ =\ e^{\frac{\pi i}{2\kappa}l^2}\
\Delta_{l,\kappa}(\lambda,\tau,\eta)\ .
\eean
\item[(ii)]
If \, $\mathrm{Im}\,\eta<0$,\, $\mathrm{Im}\,\tau>0$,\,
 $\mathrm{Im}\,\eta/\tau < 0$, then
for $l=2, \dots , \kappa-2$ we have
\bean\label{2mp}
C^-(\tau,\eta)\ e^{-\,\frac{\pi i}{2\tau}(\kappa+2)\lambda^2}
 \Delta_{l,\kappa}\left(
\frac\lambda\tau,-\,\frac1\tau,\frac\eta\tau\right)
\ =\
 \sum_{j=2}^{\kappa-2}\ 
\Delta_{j,\kappa}(\lambda,\tau,\eta)\ S^-_{j,l}(\tau,\eta)\ ,
\eean
where
\bea
S^-_{j,l}(\tau,\eta) &=&  Q\left(\frac l\kappa,\frac{\tau}{2\eta\kappa},
-\frac1{2\kappa}\right)
\ \times
\\
& & \ {} \ {} \ {} 
\left(
u\left(\frac j\kappa,-\frac l\kappa,
\frac1{2\eta\kappa},\frac\tau{2\eta\kappa},-\,\frac1{2\kappa}\right)
-u\left(\frac j\kappa, \frac l\kappa,
\frac1{2\eta\kappa},
\frac\tau{2\eta\kappa},-\,
\frac1{2\kappa}\right)\right)\ ,
\\
C^-(\tau,\eta) & = & - 2\pi i \, \sqrt\frac{2\kappa i}{\tau}\
e^{\pi i\frac{2\eta}\kappa +\frac{\pi i}\tau 4\eta^2(1 - \frac{1}{2\eta\kappa})
+\frac{\pi i}{6\kappa\tau}\psi(\tau, -2\eta\kappa)}\ ,
\\
\psi(\tau,p) &=& 2\ (\,8\eta^2 + \tau^2 + p^2 - 3p + 3\tau + 3\tau p + 1\,)\ .
\eea
Here the square root is with positive real part.
\item[(iii)]
If \, $\mathrm{Im}\,\eta < 0$, \,$\mathrm{Im}\,\tau > 0$,\,
 $\mathrm{Im}\,\eta/\tau > 0$, then
for $l=2, \dots , \kappa-2$,
\bean\label{3mp}
C^+(\tau,\eta)\ e^{-\,\frac{\pi i}{2\tau}(\kappa+2)\lambda^2}
 \Delta_{l,\kappa}\left(
\frac\lambda\tau,-\,\frac1\tau,-\frac\eta\tau\right)
\ =\
 \sum_{j=2}^{\kappa-2}\ 
\Delta_{j,\kappa}(\lambda,\tau,\eta)\ S^+_{j,l}(\tau,\eta)\ ,
\eean
where
\bea
S^+_{j,l}(\tau,\eta) &=&  Q\left(\frac l\kappa,-\frac{\tau}{2\eta\kappa},
\frac1{2\kappa}\right)
\ \times
\\
& & \ {} \ {} \ {} 
\left(
u\left(\frac j\kappa,-\frac l\kappa,
\frac1{2\eta\kappa},-\frac\tau{2\eta\kappa},\frac1{2\kappa}\right)
-u\left(\frac j\kappa, \frac l\kappa,
\frac1{2\eta\kappa},
-\frac\tau{2\eta\kappa},\frac1{2\kappa}\right)\right)\ ,
\\
C^+(\tau,\eta) & = & - 2\pi i \, \sqrt\frac{2\kappa i}{\tau}\
e^{\pi i\frac{2\eta}\kappa +\frac{\pi i}\tau 4\eta^2(1 + \frac{1}{2\eta\kappa})
+ \frac{\pi i}{6\kappa\tau}\psi(-\tau, +2\eta\kappa)}\ .
\eea
Here the square root is with positive real part.
\end{enumerate}
\end{thm}

Notice that in \Ref{2mp} and \Ref{3mp} the matrix $S$ has the property
$S(\tau - 2\eta\kappa, \eta) = S(\tau, \eta)$. Hence the right hand sides 
of formulas \Ref{1mp} - \Ref{3mp} are solutions of the qKZB equation with parameters
$\tau, \eta$. Therefore formulas \Ref{1mp} - \Ref{3mp} give transformations
of theta function solutions of the qKZB equation.


\begin{proof} Property (i) is a direct corollary of definitions. We prove (ii).
The proof of (iii) is analogous.

The proof of (ii) is based on a discrete version of the following identity
 for the function $u$ \cite{FV3}. If Im\,$(\eta\tau/p) < 0$
and Im\, $(p/\tau) > 0$, then  
\bean\label{continuous}
{}
\\
\int u(\lambda, \mu, \tau, p, \eta)\,
u\left(-\frac\mu p, \frac\nu p, -\frac 1 p, -\frac\tau p, \frac\eta p\right)\,
Q(\mu, p, \eta)\, \rho(\mu, p, \eta)\
e^{- \frac{\pi i \tau}{4\eta p}\mu^2}\, d\mu \ =
\phantom{a a a a a a a a a a }
\notag
\\
2 \pi i \, \sqrt{ \frac{4\eta \tau}{ip}}\,
\rho(\lambda, \tau, \eta)\,
\rho\left(\frac\nu  p , -\frac\tau p, \frac\eta p \right)\,
e^{\frac{\pi i p}{4\eta \tau}\,(\lambda^2 + \nu^2/p^2)}\,
e^{-\frac {\pi i \eta}{3\tau p}\,\psi(\tau, p)}\,
u\left(\frac\lambda\tau, \frac\nu\tau,-\frac1\tau, \frac p\tau, \frac\eta\tau\right)\,.
\notag
\eean
Here 
$\rho(\lambda, \tau, \eta) = e^{-\frac{\pi i}\tau\,(\lambda^2-4\eta^2)}$.
 The integration over $\mu$ is over the path 
$x \mapsto x\eta + \epsilon, \, x \in \R$, for any generic real $\epsilon$.

The discrete version has the following form. For any  admissible $l$
and generic $\epsilon$ we have
\bean\label{mod 1}
&&
{}
\\
&&\sum_{\mu\in 2\eta(\Z+\epsilon)} 
u(\lambda, \mu, \tau, -2\eta\kappa, \eta)
Q(\mu, -2\eta\kappa, \eta)
 e^{\pi i \,\frac{\tau + 4\eta}{2\kappa}\,\left(\frac\mu{2\eta}\right)^2}
Q \left(\frac{ l}\kappa,\frac{\tau}{2\eta\kappa},- \frac{1 }{2\kappa}\right)
 \times
\notag
\\
&&
\ {}\ 
u\left(\frac\mu {2\eta\kappa},- \frac l {\kappa}, \frac 1 {2\eta\kappa},
 \frac\tau {2\eta\kappa}, -\frac 1 {2\kappa}\right)
\ =\  C^-(\tau,\eta) \, e^{-\frac{\pi i}{2\tau}(\kappa+2)\lambda^2 }\,
e^{-\frac{\pi i l\epsilon}{\kappa}}\
\widehat\Delta_{l,\kappa}
\left(\frac{\lambda}{\tau}, - \frac\epsilon\kappa, -\frac{1}\tau, \frac\eta\tau
\right)\ ,
\notag
\eean
where
\bea
&&
\widehat\Delta_{l,\kappa}(\lambda,\epsilon,\tau,\eta)\ =\
e^{\frac{2\pi i\eta }\kappa\, l^2}
Q(2\eta l,-2\eta\kappa,\eta)\ \times
\ {} \ {}\ {}\   \ {} \ {}\ {}\   \ {} \ {}\ {}\   
\ {} \ {}\ {}\   \ {} \ {}\ {}\   \ {} \ {}\ {}\   
\\
&&
\ {} \ {}\ {}\   
\int_\gamma
\Omega_{2\eta}(t,\tau,-2\eta\kappa)
\frac{\theta(\lambda+t,\tau)\theta(2\eta l+t,-2\eta\kappa)}
     {\theta(t-2\eta,\tau)\theta(t-2\eta,-2\eta\kappa)}
e^{-2\pi ilt/\kappa}
\theta_{l,\kappa}\left({\textstyle\frac2\kappa}\,t-\lambda + \epsilon,\tau
\right)dt.
\eea

The proof of this identity is the same as the proof in \cite{FV3} of identity
\Ref{continuous} but instead of evaluating the Gaussian integral one
applies the Poisson summation formula. 

Identity \Ref{mod 1} implies
\bean\label{mod 2}
&&
{}
\\
&&\sum_{\mu\in 2\eta(\Z+\epsilon)} 
u(\lambda, \mu, \tau, -2\eta\kappa, \eta)
Q(\mu, -2\eta\kappa, \eta)
 e^{\pi i \,\frac{\tau + 4\eta}{2\kappa}\,\left(\frac\mu{2\eta}\right)^2}
Q \left(\frac{ l}\kappa,\frac{\tau}{2\eta\kappa},- \frac{1 }{2\kappa}\right)
 \times
\notag
\\
&&
\ {}\ \ {}\ 
[u\left(\frac\mu {2\eta\kappa},- \frac l {\kappa}, \frac 1 {2\eta\kappa},
 \frac\tau {2\eta\kappa}, -\frac 1 {2\kappa}\right)
-
u\left(-\frac\mu {2\eta\kappa},- \frac l {\kappa}, \frac 1 {2\eta\kappa},
 \frac\tau {2\eta\kappa}, -\frac 1 {2\kappa}\right)]\ =
\notag
\\
&& 
\ {}\  \ \ {}\ 
 C^-(\tau,\eta)  \, e^{-\frac{\pi i}{2\tau}(\kappa+2)\lambda^2 }\,
[e^{-\frac{\pi i l\epsilon}{\kappa}}\
\widehat\Delta_{l,\kappa}
\left(\frac{\lambda}{\tau}, - \frac\epsilon\kappa, -\frac{1}\tau, \frac\eta\tau
\right) -
e^{\frac{\pi i l\epsilon}{\kappa}}\
\widehat\Delta_{-l,\kappa}
\left(\frac{\lambda}{\tau}, - \frac\epsilon\kappa, -\frac{1}\tau, \frac\eta\tau
\right)]\ .
\notag
\eean
Let $\epsilon$ tend to zero and $\mu = 2\eta (m+\epsilon)$,\, $m\in\Z$.
Then the limit of each factor 
of the corresponding term in the left hand side of \Ref{mod 2} is well defined
unless $m \equiv \pm1\mod \kappa$. Moreover, if $m \equiv 0 \mod \kappa$, then
the limit is zero, since the function
$u\left(\frac{m} {\kappa},- \frac l {\kappa}, \frac 1 {2\eta\kappa},
 \frac\tau {2\eta\kappa}, -\frac 1 {2\kappa}\right)$ is multiplied by $-1$
 under the shift  $m \mapsto m+\kappa$.

If $m \equiv \pm1\mod \kappa$ and $\epsilon$ tends to zero, then 
$Q(2\eta(m+\epsilon), -2\eta\kappa, \eta)$ tends to infinity, but the factor
\bea
[u\left(\frac{m+\epsilon} {\kappa},- \frac l {\kappa}, \frac 1 {2\eta\kappa},
 \frac\tau {2\eta\kappa}, -\frac 1 {2\kappa}\right)
-
u\left(-\frac{m+\epsilon}{\kappa},- \frac l {\kappa}, \frac 1 {2\eta\kappa},
 \frac\tau {2\eta\kappa}, -\frac 1 {2\kappa}\right)]
\eea
tends to zero by Lemma \ref{vanishing}. Hence the limit 
of the corresponding term in the left hand side of \Ref{mod 2} is well defined
in this case as well. Taking the limit $\epsilon \to 0$ we observe that for any $r\in\Z$
the terms corresponding to $m = r\kappa +1$ and $m = r\kappa - 1$ are canceled by
Lemma \ref{vanishing} applied to the factor 
$u(\lambda, 2\eta(m+\epsilon), \tau, -2\eta\kappa, \eta)$. 
Now taking the limit $\epsilon \to 0$ in \Ref{mod 2} we get the formula of part (ii) of the theorem.
\end{proof}

\bigskip

\begin{thm}\label{A,B}
For every admissible $l$ we have
\bea
\Delta_{l,\kappa}(\lambda+1,\tau,\eta) &=&
 (-1)^{l+1} \Delta_{l,\kappa}(\lambda,\tau,\eta)\ ,
\\
e^{\pi i (\kappa + 2)(\lambda + \tau/2)}\,
\Delta_{l,\kappa}(\lambda+\tau,\tau,\eta)
&=&
\Delta_{l+\kappa,\kappa}(\lambda,\tau,\eta)\ .
\eea
\end{thm}
\bigskip

This theorem gives additional transformations of
the theta function solutions of the qKZB equation with parameters $\tau, \eta$.

\bigskip

\subsection{Modular action}
The group $\SL(2,\Z)$ is generated by matrices
\bea
T\ = \ \mdpd 1101 \ ,
\qquad
S\ =\ \mdpd 01{-1}0
\eea
satisfying the relations $S^2 = -Id,\ (ST)^3 = Id$.
 Let $\HH_+$  be the upper half plane of complex numbers
with positive imaginary part, and let $\HH_-$ be
the lower half plane. The group
$\SL(2,\Z)$ acts on $\HH_+\times \C$ with coordinates $\tau, \eta$ by the formulas
\bea
S : (\tau,\eta) \mapsto (-1/\tau, \eta/\tau) ,
\qquad
T : (\tau,\eta) \mapsto (\tau+1, \eta) .
\eea
An orbit of this action will be called admissible if it does not contain a point
$(\tau, \eta)$ with $\eta \in \R$. Let $X$ be the intersection of $\HH_+\times \HH_-$
with the union of all admissible orbits.

If $(\tau,\eta)\in \HH_+\times \C$, then the set $\{(\tau - 2\eta \kappa l),\
l = 0, 1, 2, \dots \}$ will be called the qKZB sequence generated at $(\tau,\eta)$.
If $(\tau,\eta)\in X$, then $X$ contains the qKZB sequence generated at
$(\tau, \eta)$.

Consider the space $F$ 
of all functions $u(\lambda,\tau,\eta)$, defined on $\C\times X$ such that for 
any fixed $(\tau,\eta)\in X$ we have $u \in E_\kappa(\tau,\eta)$.

Introduce transformations $A, B, T$, $S$ of functions of $F$ by
formulas
\bea
(Au)(\lambda,\tau,\eta) = u(\lambda+1,\tau,\eta) ,
&&
(Bu)(\lambda,\tau,\eta) = 
e^{\pi i (\kappa + 2)(\lambda + \tau/2)}\,
u(\lambda+\tau,\tau,\eta) ,
\\
(Tu)(\lambda,\tau,\eta) &=& u(\lambda,\tau+1,\eta) ,
\\
(Su)(\lambda,\tau,\eta) &=& C^\pm(\tau,\eta)\
e^{-\frac{\pi i}{2\tau} (\kappa + 2)\lambda^2}\
u\left(\frac{\lambda}{\tau},-\frac{1}{\tau},
\pm\frac{-\eta}{\tau}\right) ,
\eea
where in the definition of $S$ the sign $+$ is chosen if Im $\eta/\tau > 0$
and $-$ if Im $\eta/\tau < 0$.

\begin{lemma}\label{1}
 The transformations $A, B, T, S$ preserve
 the space $F$. Moreover, each of the transformations send solutions of
the qKZB equation to solutions.
\end{lemma}

The lemma is a corollary of Theorems   \ref{t-mod}     and \ref{A,B}.

\begin{lemma}[cf. \cite{FSV2}]\label{2} 
Restricted to the space $F$, the transformations $A, B, T, S$ 
satisfy the relations
\bea
A^2 = 1 ,\qquad
B^2 = 1 , \qquad S^2 = 8\pi^2\kappa ,\qquad 
(ST)^3 = 8\pi^3 \kappa^{3/2} e^{ {2\pi i}/{\kappa}-{\pi i}/{4} } ,
\qquad {}\ {}\ \\
SA = BS, \qquad AB = (-1)^\kappa BA ,\qquad
TB = - e^{\pi i\kappa/2}BAT , \qquad AT = TA .
\eea

\end{lemma}

Lemmas \ref{1}, \ref{2} in particular say that the transformations $S, T$ define 
a projective representation of $\SL(2,\Z)$ in the space of theta function 
solutions of the qKZB equation.

\bigskip

\bigskip

\section{Limiting cases}
\subsection{Conformal blocks}
It was shown in \cite{FV2}
that in
 the limit $\eta\to 0$, the qKZB heat equation \eqref{e-qKZB}
degenerates to the KZB heat equation of \cite{B}.
More precisely, if $v(\lambda,\tau,\eta)$ is a family 
of solutions
of \eqref{e-qKZB} such that $v(\lambda,\tau,\eta)=v(\lambda,\tau)
+O(\eta)$ then $v(\lambda,\tau)/\theta(\lambda,\tau)$ obeys the
KZB heat equation
with the three dimensional
irreducible representation of $sl_2$ sitting at
one marked point:
\bea\label{KZB}
2\pi i\kappa\frac
{\partial v}{\partial\tau}=
\frac
{\partial^2v}
{\partial\lambda^2}
+2\rho'(\lambda,\tau)v
+c(\tau)v, \qquad \rho=\theta'/\theta ,
\eea
for some $c(\tau)$ independent of $\lambda$.

Conformal blocks are solutions of the KZB heat equation with $c(\tau) = 0$ 
taking
values in theta functions of level $\kappa-2$.
The integral representation of conformal blocks has the following
form
\cite{FV1}:
\begin{equation}\label{e-ccb}
v_{l,\kappa}(\lambda,\tau)
=\int_0^1
\left(
\frac{\theta(t,\tau)}
{\theta'(0,\tau)}\right)^{-\frac2\kappa}
\frac{\theta(t-\lambda,\tau)\theta'(0,\tau)}
     {\theta(t,\tau)\theta(\lambda,\tau)}
\theta_{l,\kappa}\left({\textstyle{\frac2\kappa}}\, t+\lambda,\tau\right)dt
-(\lambda\to-\lambda).
\end{equation}

In the limit $\eta \to 0$ the integral representation for
elliptic hypergeometric functions turns into the integral representation for conformal blocks.
To make the statement precise we first discuss the integral representations.

The integral in \Ref{e-ccb}
is understood in the sense of analytic continuation from the range of parameters
where the exponent $-2/\kappa$ is replaced by a positive number.
We give an explicit formula for the regularization of the integral
\bea
F(\lambda,\tau)\ =\ \int_0^1 (\theta(t,\tau)/\theta'(0,\tau))^{-1-\frac 2\kappa} f(t,\lambda,\tau) dt
\eea
 where
 $f(t,\lambda,\tau) = 
\theta(t-\lambda,\tau)
\theta_{l,\kappa}\left({\textstyle{\frac2\kappa}}\, t+\lambda,\tau\right)$.

Set $\tilde f(t,\tau) = e^{ 2\pi i l t/\kappa}\theta(t,\tau)$. We have
\bea
f(t+1,\lambda,\tau) &=& - e^{2\pi i l/\kappa} f(t,\lambda,\tau),
\qquad
\tilde f(t+1,\tau)\ =\ - e^{2\pi i l/\kappa} \tilde f(t,\tau),
\\
(\tilde f(t,\tau)\theta(t,\tau)^{-1-2/\kappa})' & =&
- \frac{2}\kappa \ 
e^{ 2\pi i l t/\kappa}\,\theta(t,\tau)^{-1-2/\kappa}\,
(\theta'(t,\tau)\ -\ \pi i l\  \theta(t,\tau))\  .
\eea
The regularization of the integral can be defined as
\bean\label{regularization}
&&
F(\lambda,\tau) \ =\ 
\\
&&
\ \ \
\int_0^1 \left(\frac{\theta(t,\tau)}{\theta'(0,\tau)}\right)^{-1-\frac 2\kappa} 
( f(t,\lambda,\tau)
- \frac{ f(0,\lambda,\tau)}{\theta'(0,\tau)}\,
e^{\pi i t + 2\pi i l t/\kappa}
( \theta(t,\tau)'-\pi i l \theta(t,\tau))\, dt \ .
\notag
\eean
The added terms form a complete differential with respect to $t$ and the second factor 
is equal to zero at $t = 0, 1$, thus the integral is well defined.

The integral representation for hypergeometric theta functions has the following
form,
\bean\label{int repn q}
\Delta_{l,\kappa}(\lambda,\tau,\eta)\ =\
e^{\frac{2\pi i\eta }\kappa\, l^2} 
(I_{l,\kappa}(\lambda,\tau,\eta)-I_{l,\kappa}(-\lambda,\tau,\eta))
Q(2\eta l,-2\eta\kappa,\eta)\ ,
\eean
where 
\bean\label{start}
I_{l,\kappa}(\lambda,\tau,\eta) &=& \int_\gamma\
G(t,\tau,\eta)\ g(t,\lambda,\tau, \eta)\ dt \ ,
\\
G(t,\tau,\eta) &=&
\Omega_{2\eta}(t,\tau,-2\eta\kappa)
\frac{\theta(2\eta l+t,-2\eta\kappa)}
     {\theta(t-2\eta,\tau)\theta(t-2\eta,-2\eta\kappa)}
e^{-2\pi ilt/\kappa},
\notag
\\
g(t,\lambda,\tau, \eta) &=&
\theta(\lambda+t,\tau)
\theta_{l,\kappa}\left({\textstyle\frac2\kappa}\,t-\lambda,\tau
\right) .
\notag
\eean
Set $\tilde g(t,\tau,\eta) = e^{ 2\pi i l t/\kappa}\theta(t-2\eta,\tau)$. We have
\bea
g(t+1,\lambda,\tau,\eta) = e^{2\pi i l/\kappa} g(t,\lambda,\tau,\eta),
\qquad
\tilde g(t+1,\tau,\eta) = e^{2\pi i l/\kappa} \tilde g(t,\tau,\eta),
\eea
\bea
&&
G(t-2\eta\kappa,\tau,\eta)
\tilde g(t-2\eta\kappa,\tau,\eta)- G(t,\tau,\eta)\tilde g(t,\tau,\eta) =
\\
&&
\ {}\ \ {}\ \ {}\ \ {}\ \ {}\ \ {}\ 
G(t,\tau,\eta) e^{ 2\pi i l t/\kappa}
(e^{-4\pi i l\eta}\theta(t+2\eta,\tau) - \theta(t-2\eta,\tau))\ .
\eea
We have
\bean\label{regular}
&&
I_{l,\kappa}(\lambda,\tau,\eta) = \int_\gamma\
G(t,\tau,\eta)\ ( g(t,\lambda,\tau, \eta)-
\\
&&
\ {}\ \ {}\ \ {}\ \ {}\ \ {}\ 
\frac{g(2\eta,\lambda,\tau, \eta)}
{e^{4\pi i \eta(1-l/\kappa)}\theta(4\eta)}
 e^{ 2\pi i l t/\kappa}
(e^{-4\pi i l\eta}\theta(t+2\eta,\tau) - \theta(t-2\eta,\tau))\ dt\ .
\notag
\eean
The integrals \Ref{start} and \Ref{regular} are equal since the added terms form a discrete
differential, see \cite{FTV, FV3}.

The second factor of the integrand in \Ref{regular}
is zero at $t=2\eta$. The integrand is holomorphic 
at $t=2\eta$ since $t=2\eta$ is a simple pole of $G(t,\tau,\eta)$. By Stokes' theorem 
the integration contour $\gamma$ in \Ref{regular}
 can be replaced by the new contour
$\tilde \gamma : [0,1] \to \C,\, s \mapsto s+4\eta$.

\bigskip

\begin{thm} 
If $\Delta_{l,\kappa}$ are given by \Ref{int repn q} in which formula
\Ref{regular} is used, then
\[
\lim_{\eta\to 0}\frac{2\eta\Delta_{l,\kappa}(\lambda,\tau,\eta)}
{\theta(\lambda,\tau)}=
-
\frac{(2\pi)^{-\frac2\kappa}\,
e^{\pi i\frac{l+2}\kappa}\sin{\displaystyle\frac{2\pi}\kappa}
}
{2\sin{\displaystyle\frac{\pi(l+1)}{\kappa}}
\sin{\displaystyle\frac{\pi(l-1)}{\kappa}}}
\,\prod_{j=1}^\infty(1-e^{2\pi ij\tau})^{-3-4/\kappa}
v_{l,\kappa}(\lambda,\tau),
\]
where $v_{l,\kappa}(\lambda,\tau)$ is defined in \Ref{e-ccb}
and the integral in \Ref{e-ccb} is regularized as in \Ref{regularization}.
\end{thm}

\bigskip

The proof follows from the formula in  \cite{FV4} for asymptotics of 
 $\Omega_{2\eta}$ as $\eta$ tends to zero.

Under the limit $\eta\to 0$ the discussed in Theorem \ref{t-mod}
the modular properties of hypergeometric theta functions
$\Delta_{l,\kappa}$ 
degenerate to the calculated in \cite{FSV1} modular properties of 
conformal blocks  $v_{l,\kappa}(\lambda,\tau)$.

\bigskip

\subsection{Elliptic Macdonald polynomials}
The $A_1$-Macdonald polynomials $P^{(m)}_j(x)\in\C[x,x^{-1}]$ 
(or Askey--Wilson polynomials) are
Laurent polynomials depending on two non-negative integers $m$ and $j$ and a parameter $q$, see \cite{M}. 
They are defined by
the conditions:
\begin{enumerate}
\item[(i)] $P^{(m)}_j(x) = x^j+$ terms of lower degree, 
\item[(ii)]  $P^{(m)}_j(x^{-1}) = P^{(m)}_j(x)$,
\item[(iii)] 
$\langle P^{(m)}_j, P^{(m)}_k \rangle = 0$ if $j\neq k$,
 with respect to
the inner product $\langle P,Q \rangle =$
constant term of 
\bea
 P(x) Q(x^{-1})\prod_{j=0}^{m-1}(1 - q^{2j}x^2)(1 - q^{2j}x^{-2})\ .
\eea
\end{enumerate}

The $A_1$- (and more generally $A_n$-) Macdonald polynomials
are traces of intertwining operators for the quantum
group $U_q(sl(2))$ (resp.\ $U_q(sl(n+1))$)
\cite{EK1}, cf. \cite{EV1, EV2}. The simplest non-trivial
case is the case of $m = 2$. The Macdonald polynomials
have the form
\begin{gather}\label{e-M}
P_j^{(2)}(x,q)\ =\ \frac{x^{j+3}-a x^{j+1}+a x^{-j-1}-x^{-j-3}}
            {\Pi(x,q)},\\
\Pi(x,q)={(qx-(qx)^{-1})(x-x^{-1})(q^{-1}x-qx^{-1})},
\qquad a=\frac{q^{j+3}-q^{-j-3}}{q^{j+1}-q^{-j-1}}.\notag
\end{gather}
Let us define elliptic Macdonald polynomials for $m = 2$ by
\begin{gather}\label{e-E}
P_{j,\kappa}(x,q,p)=
e^{-{\pi i}\frac{4\eta+\tau}{2\kappa}(j+2)^2+\pi i3\tau/4}
\frac{
\Delta_{j+2,\kappa}(\lambda,\tau,\eta)}
{\theta(\lambda-2\eta,\tau)\theta(\lambda,\tau)\theta(\lambda+2\eta,\tau)},\\ 
x=e^{\pi i\lambda},\quad q=e^{-2\pi i\eta},\quad p=e^{2\pi i\tau},
\quad j=0,\dots,\kappa-4.\notag
\end{gather}
The exponential function in the definition of $P_{j,\kappa}$
ensures that the result is 1-periodic in $\tau,\eta$ and can
thus be written as a function of $p$ and $q$.

\bigskip

\begin{thm}\ %
\begin{enumerate}
\item[(i)]
The elliptic Macdonald polynomials form a basis of
$\Theta_{\kappa-4}(\tau)^\mathrm{even}$.
\item[(ii)] As $p\to 0$, we have
$
P_{j,\kappa}(x,q,p) = A_{j,\kappa}(q)P_j^{(2)}(x,q)+O(p),
$
for some $A_{j,\kappa}(q)\neq 0$.
\end{enumerate}
\end{thm}

\bigskip
\noindent{\it Proof:}
Part (i) follows from  \eqref{e-cb} and
Theorem \ref{t-1}.
Let us prove part (ii).
As $\tau\to i\infty$ ($p\to 0$) the denominator in
\eqref{e-E} is 
$i e^{3\pi i\tau/4}
(\Pi(x,q)+O(p))$. The numerator 
$\Delta_{j+2,\kappa}(\lambda,\tau,\eta)$ behaves as
\[
e^{{\pi i}\frac{4\eta+\tau}{2\kappa}(j+2)^2}
Q(2\eta (j+2),-2\eta\kappa,\eta)
(u(\lambda,2\eta(j+2),i\infty,-2\eta\kappa,\eta)
-(\lambda\to-\lambda)+O(p)),
\]
as $p\to 0$. Here,
\bea
u(\lambda,2\eta(j+2),i\infty,-2\eta\kappa,\eta)
&=&
e^{-{\pi i\lambda(j+2)}}
\int_\gamma
\prod_{j=0}^\infty
\frac
{1-q^{2j\kappa+2}e^{2\pi it}}
{1-q^{2j\kappa-2}e^{2\pi it}}\\
&\times&
\frac{\sin\pi(\lambda+t)\theta(2\eta(j+2)+t,-2\eta\kappa)}
     {\sin\pi(t-2\eta)\theta(t-2\eta,-2\eta\kappa)}dt\ .
 \eea
Considering the right hand side of this formula as a function of
$x$, we see that $\lim_{\tau\to i\infty}P_{j,\kappa}$ 
has the form \eqref{e-M} for some coefficient $a$
up to a factor independent of $x$. 
 The fact
that the numerator is divisible by the denominator
determines the value of $a$ uniquely.
The fact that $A_j(q) \neq 0$ follows from the
limiting case of the inversion relation \eqref{e-ir}.
\hfill$\square$

\medskip

Let us now study the qKZB heat equation in this limit.
The integration kernel appearing in the equation involves
the limit of $u$ as both $\tau$ and $\sigma$ tend to $i\infty$:
\begin{eqnarray*}
u(\lambda,\mu,i\infty,i\infty,\eta)
&=&
e^{{\textstyle-\frac{i\pi\lambda\mu}{2\eta}}}
\int_\gamma
\frac{\sin\pi(\lambda+t)\sin\pi(\mu+t)}
     {\sin\pi(t+2\eta)\sin\pi(t-2\eta)}dt\\
&=&
i\,e^{{\textstyle-\frac{i\pi\lambda\mu}{2\eta}}}
\frac{(q^{-2}\cos\pi(\lambda+\mu)-\cos\pi(\lambda-\mu))}
{\sin{4\pi\eta}}\ .
\notag
  \end{eqnarray*}
Combining this with
\[
Q(\mu,i\infty,\eta)=\frac{\pi\sin4\pi\eta}
{\sin\,\pi(\mu-2\eta)\,\sin\,\pi(\mu+2\eta)},
\]
we obtain in the limit the Macdonald-Mehta type identity
\bean\label{mehta}
q^{-\frac{(j+2)^2}2}\,P^{(2)}_{j}(e^{\pi i\lambda},q)\ =\
\int_{\eta\R}\ V(\lambda,\mu)\ P^{(2)}_{j}(e^{-i\pi\mu},q)\ d\mu\ ,
\eean
\[V(\lambda,\mu)\ =\ 
\frac{q^{-2}e^{{\textstyle-\frac{i\pi(\lambda+\mu)^2}{4\eta}}}}
{2i\sqrt{4i\eta}}\ 
\frac{(q^{-2}\cos\pi(\lambda+\mu)-\cos\pi(\lambda-\mu))
\,\sin\pi\mu}
{\sin\pi(\lambda-2\eta)\,\sin\pi\lambda
\,\sin\pi(\lambda+2\eta)}\ ,
\]
cf. \cite{EV2}.

The limit of the difference equation \Ref{difference} has the following form.
Define the infinite order difference operator
\bean\label{operator}
&&
T(q) = 
\\
&& \ {} \ {} \ {} \ {} 
-q^{-2} \sum_{m\in\Z}
\frac{q^{-2}\cos\,\pi(2\eta m)-\cos\,\pi(2\lambda-2\eta m)}
{\sin\,\pi(\lambda-2\eta)\,\sin\,\pi\lambda
\,\sin\pi(\lambda+2\eta)}\ \sin\,\pi(\lambda-2\eta m)\,
q^{-\frac{m^2}2}\,T_{-2\eta m}
\notag
\eean
where $T_{-2\eta m} : v(\lambda) \mapsto 
v(\lambda -2\eta m)$ is the operator of the shift of the argument by
$-2\eta m$. Set 
\bea
\theta_0(x,q)\ =\ \sum_{m\in\Z} x^m q^{-\frac{m^2}2}\ .
\eea
Then
\bean\label{limit of diff eqn}
T(q)\ P^{(2)}_{j}(e^{\pi i\lambda},q)\ = \ 
\theta_0(q^{j+2},q)\,
P^{(2)}_{j}(e^{\pi i\lambda},q)\ .
\eean


\bigskip

\subsection{Remarks on the operator $T(q)$}
Let $x = e^{\pi i \lambda}, q = e^{-2\pi i \eta}$ as before.
Following \cite{M} introduce  operators acting on Laurent
polynomials: \  $w : f(x) \mapsto f(qx^{-1})$,\,
$\Gamma : f(x) \mapsto f(qx)$,
\bea
Y \ =\
\frac{q^{-1}-qx}{1-x}\Gamma\ +\
\frac{q-q^{-1}}{1-x} w\ .
\eea
Then 
\bea
Y^{-1} \ =\
\Gamma^{-1}\frac{q-q^{-1}x}{1-x}\ +\
w x \frac{q-q^{-1}}{1-x} \ .
\eea
For any symmetric Laurent polynomial $f(x)$, $f(x)=f(x^{-1})$, and any
$j = 0, 1, \dots$,  one has
\bean\label{chered}
f(Y)\,P^{(2)}_j(x)\ =\ f(q^{j+2})\, P^{(2)}_j(x)\ ,
\eean
see \cite{M}.
The Macdonald polynomials form a basis in the space of symmetric
Laurent polynomials. Comparing \Ref{limit of diff eqn}    and \Ref{chered}
we see that
\bea
T(q)  \ =\ \theta_0(Y,q)
\eea
as operators on the space of symmetric Laurent polynomials.

\subsection{Orthogonality relation} 
In terms of elliptic Macdonald polynomials 
the inversion relation \Ref{inversion} says that
\bean\label{orthogonality}
&&
\frac1{32\pi^2\eta}\,\int_0^2
\,u(-\mu,2\eta l,\tau,-2\eta\kappa,-\eta)\,
P_{j,\kappa}(e^{\pi i\mu},q,p)\,
\theta(\mu,\tau)\,d\mu\
\ {}\ \ {}\ \ {}\ \ {}\ \ {}\ \ {}\ \ {}\ \ {}\ \ {}\ \ {}\ \ {}\ \ {}\ 
\\
&&
\ {}\ \ {}\ \ {}\ \ {}\ \ {}\ \ {}\ \ {}\ \ {}\ \ {}\ \ {}\ \ {}\ \ {}\ 
\ {}\ \ {}\ \ {}\ \ {}\ \ {}\ \ {}\ \ {}\ \ {}\ \ {}\ \ {}\ \ {}\ \ {}\ 
\ {}\ \ {}\ \ {}\ \ {}\ \ {}\ \ {}\ \ {}\ \ {}\ \ {}\ \ {}\ \ {}\ \ {}\ 
=\
\delta_{l,j}\,e^{\pi i\frac{4\eta+\tau}{2\kappa}j^2}\ 
\notag
\eean
for $l,j= 2, \dots , \kappa-2$.
Note that the integrand is a 2-periodic entire function of $\mu$. Thus the integral is equal 
to the constant Fourier
coefficient of the integrand multiplied by 2.

We want to present this integral as a pairing of suitable functional spaces. 
Namely, 
consider the vector space $F$ of entire 2-periodic functions of $\lambda$.
Define the subspace $H$ as the subspace generated by the functions of the form
$ h(\lambda)- h(\lambda+2\tau)
e^{-2\pi i (\kappa-2)(\tau+\mu)}$ and $ h(\lambda)+h(-\lambda)$ where
$h \in F$.
Then the map 
\bea
f\otimes g \ \mapsto\ 
\int_0^2 f(\mu)\,g(\mu)\,\theta(\mu,\tau)\,d\mu\
\eea
defines a perfect pairing 
$(F/H) \otimes \Theta_{\kappa-4}(\tau)^\mathrm{even} \to \C$. Formula \Ref{orthogonality}
says that the functions $u(-\mu,2\eta l,\tau,-2\eta\kappa,-\eta)$, $l = 2, \dots , \kappa-2$, 
considered as elements of
$F/H$, form a basis dual to the basis of elliptic Macdonald polynomials up to 
multiplication by scalars.


\end{document}